\documentclass{amsart}
\usepackage{amsmath, amsthm}
\usepackage{amsmath,amscd}
\usepackage{amsfonts}
\usepackage{amssymb, latexsym}
\usepackage{url}

\newtheorem{theorem}{Theorem}

\newtheorem{lemma}[theorem]{Lemma}
\newtheorem{definition}[theorem]{Definition}

\def\gp{{\rm gp}}
\author[Delaram Kahrobaei]{Delaram Kahrobaei}
\address{Delaram Kahrobaei, Doctoral Program in Computer Science, CUNY Graduate Center,
365 Fifth Avenue, New York, NY 10016\\
Mathematics Department, New York City College of Technology (CUNY)\\
300 Jay Street, Brooklyn, NY 11201}
\email{dkahrobaei@gc.cuny.edu\\ https://wfs.gc.cuny.edu/DKahrobaei/www/}

\author[Stephen Majewicz]{Stephen Majewicz}
\address{Stephen Majewicz, Mathematics Department, Kingsborough Community College
(CUNY)\\
2001 Oriental Blvd., Brooklyn, NY 11235}
\email{smajewicz@kbcc.cuny.edu}

\thanks{The research of the first author has been supported by the PSC CUNY grant from the research foundation of the City University of New York and the City Tech Foundation.}
\begin{document}

\title[Residual Solvability]{On the residual solvability of
generalized free products of solvable groups}

{\abstract{In this paper, we study the residual solvability of the
generalized free product of solvable groups.}}

\maketitle
\centerline{Dedicated to Laci Babai}

\section{Introduction and Motivation}

In his celebrated paper \cite{PH54}, Hall introduced the residual
properties. Let $\mathcal{P}$ be a property of groups, and let $G$
be a group. We say that $G$ is \emph{residually $\mathcal{P}$} if
for any non-identity element $g \in G,$ there exists a
homomorphism $\psi$ from $G$ to a group having property
$\mathcal{P}$ such that $\psi(g) \neq 1.$

We recall the definition of the generalized free products of groups. Let there be given a nonempty set of groups $\{G_\lambda | \lambda \in \Lambda\}$ together with a group $H$ which is isomorphic with a subgroup $H_\lambda$ of $G_\lambda$ by means of monomorphism $\phi_\lambda: H \rightarrow G_\lambda  (\lambda \in \Lambda)$. There is an exceedingly useful object known as the free product of the $G_\lambda$'s with the amalgamated subgroup $H$. Roughly speaking this is the largest group generated by the $G_\lambda$'s in which the subgroups $H_\lambda$ are identified by means of the $\phi_\lambda$. Generically such groups are known as generalized free products (see \cite{DR}).

In the literature, residual finiteness has attracted much
attention. In \cite{GB63}, for example, Baumslag studies the
residual finiteness of generalized free products of nilpotent
groups. Another study pursued by M. Sapir \cite{S10} and D. Wise
\cite{DW02} focuses on the residual finiteness of one-relator
groups.

On the other hand, the exploration of residually solvable groups
is in its early stages. One classical work on the residual
solvability of positive one-relator groups is due to Baumslag (see
\cite{GB71}). The first author studied doubles of residually
solvable groups (see \cite{DK2005}), as well as the residual
solvability of the generalized free products of finitely generated
nilpotent groups (see \cite{DK2010}). In \cite{ADK}, Arzhantseva,
de la Harpe, and the first author introduce the true prosolvable
completion of a group and give several interesting open problems
that motivate the importance of residually solvable groups.

In this paper, we generalize some of the results in \cite{DK2010}
which carry over when nilpotent groups are replaced by solvable
groups.

\section{Results}

In this section, we state our main results. Their proofs will be
presented in Section $4.$

\vspace{.1in}

We begin by proving that the generalized free product of a
nilpotent group and a solvable group is not necessarily perfect.

\begin{theorem} \label{Not_Perfect}
The generalized free product of a nilpotent group and a solvable
group amalgamated by a proper subgroup of them is not perfect.
\end{theorem}

Next we consider the case where the amalgamated subgroup of the
generalized free product of two solvable groups is cyclic. By
choosing an appropriate solvable filtration of each factor so that
the generator of the amalgamated subgroup does not lie in the
$n^{th}$ term of the derived series, it follows that the resulting
amalgam is residually solvable.

\begin{theorem} \label{Cyclic_amalgam}
The generalized free product of two solvable groups amalgamating a
cyclic subgroup is residually solvable.
\end{theorem}

In \cite{DK2010}, Kahrobaei proved the generalized free product of
two finitely generated solvable groups amalgamated by central
subgroups is a solvable extension of a residually solvable group
and, hence, is residually solvable. We generalize this result for
any finite number of solvable groups.

\begin{theorem} \label{finiteNilpotentRS}
The generalized free product of a finite number of solvable
groups, amalgamating a central subgroup in each of the factors, is
residually solvable.
\end{theorem}

Our next theorem deals with the generalization of doubles (recall
that a \emph{double} is the amalgamated product of two groups
whose factors are isomorphic and the amalgamated subgroups are
identified under the same isomorphism).

\begin{theorem}\label{finite_double}
Let $\{A_i \, | \, i \in I \}$ be an arbitrary indexed family of
isomorphic solvable groups such that $\bigcap_{i \in I} A_i = C,$
and let $G$ be the generalized free product of the $A_i$'s
amalgamated by $C.$ Then $G$ is residually solvable.
\end{theorem}

Our final result pertains to generalized free products containing
a finitely generated torsion-free abelian factor.

\begin{theorem} \label{abelianonefactor}
The generalized free product of a finitely generated torsion-free
abelian group and a solvable group is residually solvable.
\end{theorem}

Note that Theorem~\ref{abelianonefactor} is a special case where
the amalgamated subgroup is central in only one of the factors.

\section{Background and Preliminary Results}

In this section, we present some preliminary material which will
be needed for the proofs of our results. We begin with a theorem
due to Neumann \cite{HN49} which plays a crucial role in our work.

\begin{theorem}[Neumann]\label{Neumann}
If $K$ is a subgroup of $G =\{ A \ast B \, ; \, C \}$, then $K$ is
an HNN-extension of a tree product in which the vertex groups are
conjugates of subgroups of either $A$ or $B$ and the edge groups
are conjugates of subgroups of $C.$ The associated subgroups
involved in the HNN-extension are also conjugates of subgroups of
$C.$ If $K \cap A = \{1\} = K \cap B,$ then $K$ is free. If $K
\cap C = \{1\},$ then $K = {\prod_{i \in I}}^* X_i \ast F,$ where
the $X_i$ are conjugates of subgroups of $A$ and $B$ and $F$ is
free.
\end{theorem}

The next result is used in the proof of Theorem~\ref{Not_Perfect}.
If $G$ is any group, let $G_{\rm ab}$ denote the abelianization of
$G.$

\begin{lemma}[Kahrobaei \cite{DK2010}] \label{onto_abelianization}
Let $G$ be the amalgamated product of $A$ and $B$ with $C_A$ and
$C_B$ identified, $$G = \{A \ast B \, ; \, C_A = C_B\}.$$ Then
$G_{\rm ab}$ maps onto
\begin{eqnarray*}
D = A_{\rm ab}/{{\gp({\bar{c}}_a \, | \, c_a \in C_A)}} \times
B_{\rm ab}/{\gp({{\bar{c}}_b}^{\, -1} \, | \, c_b \in C_B)},
\end{eqnarray*}
where $\overline{c}_a$ and $\overline{c}_b$ are the images of
$c_a$ and $c_b$ in $A_{\rm ab}$ and $B_{\rm ab},$ respectively.
\end{lemma}

Theorem~\ref{Not_Perfect} also relies on the next well-known
theorem about the Frattini subgroup $\Phi(G)$ of a nilpotent group
$G.$ Let $\delta_i G$ denote the $i^{th}$ derived subgroup of a
group $G.$

\begin{theorem}[Hirsch]\label{Frattini}
If $G$ is a nilpotent group, then $\delta_2 G \leq \Phi(G).$
\end{theorem}

The proof of Theorem 3 makes use of the generalized central
product of a set of groups. Let $Z(G)$ denote the center of a
group $G.$

\begin{definition} \label{gen_cent_prod_def}
Suppose that
\begin{eqnarray*}
\{ A_i = \langle X_i \, ; \, R_i \rangle \, | \, i \in I\}
\end{eqnarray*}
is an indexed family of presentations, and let $C$ be a group
equipped with monomorphisms
\begin{eqnarray*}
\phi_i : C \rightarrow A_i \; \text{ and } \; C \leq Z(A_i) \, \,
(\text{for all } i \in I).
\end{eqnarray*}
The group $A$ defined by the presentation
\begin{eqnarray*}
A = \langle  \cup X_i \, ; \, \cup R_i \cup \{ c \phi_i c^{-1}
\phi_j \, | \, c \in C, \, i, \, j \in I \} \cup \{ [x_i, \, x_j]
= 1, \; i, \, j \in I \} \rangle,
\end{eqnarray*}
where we assume that the $X_i$ are disjoint, is termed the
\emph{generalized central product} of the $A_i$ amalgamating the
central subgroup $C.$ Thus,
\begin{eqnarray*}
A = {\prod_{i \in I}}^{\times} \{A_i \, ; \, C\}.
\end{eqnarray*}
\end{definition}

It is not hard to show that
\begin{eqnarray*}
A = { \left( {\prod_{i \in I}}^{\times} A_i \right)}/{gp(c \phi_i
c^{-1} \phi_j \, | \, i, \,j \in I, \; c \in C) }.
\end{eqnarray*}

By a theorem of von Dyck, there are canonical homomorphisms $\mu_i
: A_i \rightarrow A.$

\begin{lemma} \label{central_any_factor}
Each homomorphism $\mu_i$ is a monomorphism.
\end{lemma}

We note that the generalized central product of finitely many
solvable groups is solvable. For the case of abelian factors, the
generalized central product of an arbitrary number of abelian
groups is abelian.

\medskip

A key ingredient for proving Theorem 6 is the following:

\begin{lemma} [Kahrobaei \cite{DK2010}] \label{double_key_lemma}
Suppose $\phi : A \rightarrow B$ is an isomorphism between two
groups $A$ and $B.$ Let $C < A,$ and let $D$ be the amalgamated
product of $A$ and $B$ amalgamating $C$ with $C \phi:$
\begin{eqnarray*}
D = \{ A \ast B \, ; \, C = C \phi \}
\end{eqnarray*}
There exists a homomorphism $\psi$ from $D$ onto one of the
factors with kernel
\begin{eqnarray*}
\gp(a (a \phi)^{-1} \, | \, a \in A),
\end{eqnarray*}
and $\psi$ injects into each factor.
\end{lemma}

\section{Proofs of Main Results}

\subsection{Proof of Theorem \ref{Not_Perfect}}

Let $A$ be a nilpotent group, and let $B$ be a solvable group.
Suppose $C_A$ and $C_B$ are proper subgroups of $A$ and $B,$
respectively. We want to show
that $G = \{ A \ast B \, ; \, C_A = C_B \}$ is not perfect.\\
\indent For any $c_a \in C_A$ and $c_b \in C_B,$ let
$\overline{c}_a$ and $\overline{c}_b$ denote their images in
$A_{\rm ab}$ and $B_{\rm ab},$ respectively. Set $\overline{C}_A =
\gp(\overline{c}_a \, | \, c_a \in C_A)$ and $\overline{C}_B =
\gp({\overline{c}_b}^{\, -1} \, | \, c_b \in C_B).$ By
Lemma~\ref{onto_abelianization}, $G_{\rm ab}$ maps onto
\begin{eqnarray*}
D = A_{\rm ab}/\overline{C}_A \times B_{\rm ab}/\overline{C}_B.
\end{eqnarray*}
We claim that $A_{\rm ab}/\overline{C}_A \neq \{1\}.$ By
Theorem~\ref{Frattini}, $\gp(C_A , \, \delta_2 A)$ is a proper
subgroup of $A$ since $C_A$ is. Thus,
$$A_{\rm ab}/\overline{C}_A \simeq (A/\delta_2 A)/(\gp(C_A, \,
\delta_2 A)/\delta_2 A) \simeq A/\gp(C_A, \, \delta_2 A)
\not\simeq \{1\},$$ proving our claim. It follows that $D \neq
\{1\}$ and, hence, $G_{\rm ab} \neq \{1\}.$ Therefore, $G$ is not
perfect.

\subsection{Proof of Theorem \ref{Cyclic_amalgam}}

Suppose $A = gp(a_{i} \ | \ i \in I)$ and $B = gp(b_{j} \ | \ j
\in J)$ are solvable groups, and let $a \in A$ and $b \in B$ be
non-identity elements. There exists positive integers $m$ and $n$
such that $a \in \delta_{m} A \backslash \delta_{m + 1} A$ and $b
\in \delta_{n} B \backslash \delta_{n + 1} B.$ Let $G = \{A \ast B
\, ; \, a = b \},$ and let $D$ be the central product of
$A/{\delta_{m + 1} A}$ and $B/{\delta_{n + 1} B}$ amalgamating $a
\delta_{m + 1} A$ with $b \delta_{n + 1} B:$
\begin{eqnarray*}
D = \{ {A/{\delta_{m + 1} A}} \times B/{\delta_{n + 1} B} \, ; \,
a \delta_{m + 1} A = b \delta_{n + 1} B \}.
\end{eqnarray*}
Let $\phi : G \rightarrow D$ be the homomorphism which maps
$a_{i}$ to $(a_{i}\delta_{m + 1} A, \, \delta_{n + 1}B)$ ($i \in
I$) and $b_{j}$ to $(\delta_{m + 1}A, \, b_{j}\delta_{n + 1} B)$
($j \in J$), and let $K$ be the kernel of $\phi.$ Then $K \cap C =
\{1\}$, where $C = \gp(a)=\gp(b)$. By Theorem~\ref{Neumann}, $K$
is a free product of conjugates of subgroups of $A$ and $B,$ and a
free group. Thus $K$ is residually solvable and, hence, $G$ is an
extension of a residually solvable group by a solvable group.
Therefore, $G$ is residually solvable.

\subsection{Proof of Theorem \ref{finiteNilpotentRS}}

Suppose that $\{A_i \, | \, i \in I\}$ is a finite indexed family
of solvable groups, and let
\begin{eqnarray*}
G = \left\{{\prod_{i \in I}}^{\ast} A_i \, ; \, C \right\}.
\end{eqnarray*}
Let $S$ be the generalized central product of the $A_i$:
\begin{eqnarray*}
S = \left(\prod_{i \in I}^{\times} A_i \right)/{\gp(c \phi_i \;
c^{-1} \phi_j \, | \, i, \, j \in I, \, c \in C)},
\end{eqnarray*}
where the $\phi_k$'s are as in Definition~\ref{gen_cent_prod_def}.
For each $i \in I,$ let $\pi_{i}$ denote the canonical
homomorphism from $S$ to $A_{i}.$ Since the $\pi_{i}$ coincide on
$C,$ they can be naturally extended to a homomorphism $\mu$ from
$G$ into $S.$ Let $K$ be the kernel of $\mu,$ and observe that
$G/K$ is solvable since $S$ is. Now, $\mu$ is one-to-one
restricted to each factor; that is,
\begin{eqnarray*}
K \cap A_i = \{1 \} \text{ for all } i \in I.
\end{eqnarray*}
By Theorem~\ref{Neumann}, $K$ is free and, consequently, $G$ is a
solvable extension of a free group. Therefore, $G$ is residually
solvable.

\subsection{Proof of Theorem \ref{finite_double}}

Let $\{A_i \, | \, i \in I \}$ be an arbitrary indexed family of
isomorphic solvable groups such that $\bigcap_{i \in I} A_i = C,$
and let $G$ be the generalized free product of the $A_i$'s
amalgamated by $C.$ Pick an $i \in I,$ and let $\phi : G
\rightarrow A_i$ be an epimorphism with kernel $K.$ Then $K$ is
free because $\phi$ restricted to each factor is injective by
Lemma~\ref{double_key_lemma}, and
\begin{eqnarray*}
A_i \cap K = \{1 \} \;\; \hbox{ for all } i \in I.
\end{eqnarray*}
Therefore, $G$ is a solvable extension of a free group and is,
thus, residually solvable.

\subsection{Proof of Theorem \ref{abelianonefactor}}

Let $A$ be a finitely generated torsion-free abelian group, and
let $B$ be a solvable group. Suppose
\begin{eqnarray*}
G = \{A \ast B \, ; \, C\}.
\end{eqnarray*}
Since $A$ is finitely generated torsion-free abelian, $C$ is a
direct factor of a subgroup $A_1$ of $A$ of finite index, i.e.
\begin{eqnarray*}
A_1 = C \times H \;\; \text{ where } [A : A_1] = n < \infty.
\end{eqnarray*}
Let $\{a_i\}$ be a distinct set of coset representatives of $A_1$
in $A;$ that is, $A/{A_1} = \bigcup_{i = 1}^n a_i A_1.$ There
exists an epimorphism from $G$ onto $\bigcup_{i = 1}^n a_i A_1$
with kernel
\begin{eqnarray*}
K = \gp_G (B, \, A_1) = \left\{{\prod_{i = 1}^n}^{\ast} B^{a_i}
\ast A_1 \, ; \, C \right\}.
\end{eqnarray*}
Now let $D$ be the normal closure of finitely many copies of $B$
in $K$:
\begin{eqnarray*}
D = \gp_K \left( \bigcup_{i = 1}^t B^{a_i} \right) =
\left\{{{\prod^t_{i=1}}}^* B^{a_i}; C \right\}.
\end{eqnarray*}
By Theorem~\ref{finite_double}, $D$ is residually solvable.
Consequently, $K$ is an extension of a residually solvable group
by an abelian group. Therefore, $G$ is (residually
solvable)-by-abelian and, hence, residually solvable.

\bibliographystyle{amsplain}

\end{document}